# Structured Multi—Matrix Variate, Matrix Polynomial Equations: Solution Techniques

Garimella Rama Murthy,
Associate Professor,
IIIT—Hyderabad,
Gachibowli, HYDERABAD-32, AP, INDIA

**ABSTRACT**

In this research paper, structured bi-matrix variate, matrix quadratic equations are considered. Some lemmas related to determining the eigenvalues of unknown matrices are proved. Also, a method of determining the diagonalizabe unknown matrices is provided. The results are generalized to multi-matrix variate, matrix polynomial equations. Briefly generalization to tensor variate polynomial equations is discussed. It is hoped that the results lead to important contributions in "non-commutative algebra".

**1. Introduction:**

Based on the efforts of many mathematicians, algebra as a branch of mathematics with rich collection of results is established for ever. It found many applications in science, engineering, economics and other fields. Particularly linear algebra established with many formal results, procedures / algorithms etc finds many applications in engineering and science. As a natural generalization of a system of linear equations, researchers formulated the problem of existence / uniqueness / cardinality of solutions of matrix quadratic equations of the form

$$X^2 A_2 + X A_1 + A_0 \equiv \bar{0}, \quad where$$

X is an unknown matrix and $\{A_2, A_1, A_0\}$ are known coefficient matrices. Gantmacher and other mathematicians studied the solution technique of such equations using the associated "polynomial matrix" [Gan]

$$\mu^2 A_2 + \mu A_1 + A_0 \equiv \bar{0}, \quad where$$

$\mu$ is a real / complex valued scalar quantity. Mathematicians also generalized the results to arbitrary matrix polynomial equations of the following form

$$X^m A_m + X^{m-1} A_{m-1} + \cdots + X^2 A_2 + X A_1 + A_0 \equiv \bar{0}.$$

The author considered arbitrary matrix power series equations and arrived at some interesting results in [Rama1]. Infact the unique solution of a structured matrix power series equation ( arising in structured Markov chains of G/M/1-type ) was computed in closed form ( i.e. the spectral representation / Jordan Canonical form was computed ) in [Rama2].

Mathematicians contributed large body of literature in the research area of "commutative algebra". Our efforts, in this research paper is to contribute to the area of "non-commutative algebra". As a natural curiosity, the author conceived the problem of computing the solution of bi-matrix variate, matrix quadratic equations. In that effort, he realized that the results applicable to arbitrary uni-matrix variate, matrix quadratic equations donot generalize to the bi-matrix variate case directly. Thus, the solution of "structured" class of bi-matrix variate, matrix quadratic equations are considered and some interesting results are derived. The logic behind the results in the structured bi-matrix variate, matrix quadratic equation case naturally generalize to the structured multi-matrix variate, matrix polynomial equations.

We note that the problem of existence/uniqueness/cardinality of solutions of multi-scalar variate polynomial equations is the central goal of algebraic geometry. Elegant results such as the Bezout's theorem are already proved and established. It is easy to see that multi-matrix variate, matrix polynomial equations constitute "structured" multi-scalar variate polynomial equations. Thus, the results derived in this research paper could have important implications to algebraic geometry ( non-commutative algebra ).

This research paper is organized as follows. In Section 2, well know results for uni-matrix variate, matrix quadratic equations are briefly summarized. Also, some interesting results related to structured class of bi-matrix variate, matrix quadratic equations are derived. In Section 3, it is shown that the results naturally generalize to the case of structured multi-matrix variate, matrix polynomial equations. In Section 4, some ideas on solving structured mutli-tensor variate, tensor polynomial equations are briefly discussed. In Section 5, conclusions are provided.

## 2. Bi-Matrix Variate, Matrix Quadratic Equations : Solution Techniques:

In this section, we first discuss the well known results associated with uni-matrix variate, matrix quadratic equations.

Consider an arbitrary uni-variate matrix quadratic equation of the following form:
$$X^2 A_2 + X A_1 + A_0 \equiv \bar{0} , where \quad \dots\dots\dots\dots\dots\dots\dots\dots\dots\dots(2.1)$$
X is an unknown matrix and $\{ A_2, A_1, A_0 \}$ are known coefficient matrices.

### Eigenvalues of all Solutions:

It is easy to see that the following factorization of the associated polynomial matrix holds true:
$$\mu^2 A_2 + \mu A_1 + A_0 \equiv ( \mu I - X )( \mu A_2 + X A_2 + A_1 ). \quad \dots(2.2)$$
Thus, it is clear that the eigenvalues of all possible solutions of (2.1) are the zeroes of the determinental polynomial
$$Det ( \mu^2 A_2 + \mu A_1 + A_0 ) \dots\dots\dots\dots\dots\dots\dots\dots\dots\dots\dots\dots\dots\dots\dots(2.3)$$
Thus, it is evident that there are $\binom{2n}{n}$ equivalence classes ( equivalent in the sense of sharing a common set of eigenvalues ) of solutions.

### Matrix Solution in Each Equivalence Class:

Let a solution in an equivalence class be represented as ( Jordan Form )
$$X = T D T^{-1}.$$
Substituting in (2.1), we have that
$$T D^2 T^{-1} A_2 + T D T^{-1} A_1 + A_0 \equiv \bar{0} \dots\dots\dots\dots\dots\dots\dots\dots(2.4)$$
Equivalently, we have that
$$T [ D^2 T^{-1} A_2 + D T^{-1} A_1 + T^{-1} A_0 ] \equiv \bar{0} \dots(2.5)$$
Since, the matrix 'T' is non-singular, we have that
$$[ D^2 T^{-1} A_2 + D T^{-1} A_1 + T^{-1} A_0 ] \equiv \bar{0} \dots..(2.6)$$

Thus, the solutions in each class are determined as the solution of a set linear equations of the above form.

Keeping the above results in mind, we now consider bi-matrix variate, matrix quadratic equations. We quickly realize that we need to restrict consideration to a "structured" class of bi-matrix variate, matrix quadratic equations so that some reasonable results ( in the spirit of the above discussion ) can be derived.

Since, matrix multiplication is generally not commutative, the relative placement of unknown matrix variables with respect to the known coefficient matrices leads to the following interesting cases ( i.e. each bi-matrix variate monomial belongs to the following cases ):

(i) Both the unknown matrix variables ( X, Y ) occur on the left hand side of the known coefficient matrices. Also, they have common left eigenvectors corresponding to possibly different eigenvalues

(ii) Both the unknown matrix variables ( X, Y ) occur on the right hand side of the known coefficient matrices. Also, they have common right eigenvectors corresponding to possibly different eigenvalues

(iii) The unknown matrices occur on both the sides of the known coefficient matrices. Also, they have common left and right vectors corresponding to possibly different eigenvalues.

**Note:** All the other cases can in principle occur. But deriving any interesting results on the matrix solutions of such matrix quadratic equations donot seem to be possible. Thus, in this research paper we only consider the following bi-matrix variate, matrix quadratic equations:

(i) $X^2 A + Y^2 B + X Y C + X D + Y E + F \equiv \bar{0}$..........(2.7)
(ii) $A X^2 + B Y^2 + C X Y + D X + E Y + F \equiv \bar{0}$.........(2.8)
(iii) $X A X + Y B Y + X C Y + X D + Y E + F \equiv \bar{0}$.........(2.9)

Now, we take a closer look on the assumption about eigenvectors in cases (i) and (ii).

**Lemma 1:** In the cases (i), (ii) considered above, when X, Y are diagonalizable, they share common eigenvectors and the matrices commute

**Proof:** Consider the case where the unknown matrices { X, Y } are both diagonalizable. Thus, in this case, the matrix of right/left eigenvectors effectively determines the matrix of left/right eigenvectors. Hence we have that
$X = T D_1 T^{-1}$ and $Y = T D_2 T^{-1}$, where
the matrices $D_1, D_2$ contain the eigenvalues of matrices X, Y respectively. Also we have that
$$X Y = T D_1 T^{-1} T D_2 T^{-1} = T D_1 D_2 T^{-1}$$
$$= T D_2 D_1 T^{-1} = T D_2 T^{-1} T D_1 T^{-1} = Y X \ldots(2.10)$$

Thus, in this case the diagonalizable matrices X, Y commute.     Q.E.D.

**Lemma 2:** Consider the bi-matrix variate, matrix quadratic equation of the form in (i) above i.e.
$$X^2 A + Y^2 B + XYC + XD + YE + F \equiv 0$$
Suppose the matrices $\{X, Y\}$ have common left eigenvectors corresponding to all the eigenvalues $\{\alpha, \mu\}$ (not necessarily same). Then, we necessarily have that all the eigenvalues of unknown matrices $\{X, Y\}$ are zeroes (solutions) of the determinental polynomial (bi-scalar variate)

$$\text{Det}(\alpha^2 A + \mu^2 B + \alpha\mu C + \alpha D + \mu E + F) \quad \ldots\ldots\ldots\ldots(2.11)$$

**Proof:** From the hypothesis, we have that

$$\bar{g} X = \alpha \bar{g} \quad \text{and} \quad \bar{g} Y = \mu \bar{g} \quad \text{for every eigenvalue } \alpha \text{ of unknown}$$

matrix X and corresponding eigenvalue $\mu$ of unknown matrix Y.

Hence, we have that
$$\bar{g}(X^2 A + Y^2 B + XYC + XD + YE + F) \equiv \bar{0} \quad \ldots\ldots(2.12)$$
Thus,
$$\bar{g}(\alpha^2 A + \mu^2 B + \alpha\mu C + \alpha D + \mu E + F) \equiv \bar{0} \ldots\ldots(2.13)$$

Hence, we have that
$$\text{Det}(\alpha^2 A + \mu^2 B + \alpha\mu C + \alpha D + \mu E + F) = 0 \quad \ldots(2.14)$$

Thus, the eigenvalues of the unknown matrices X, Y are necessarily the zeroes of bi-scalar variate determinental polynomial:

$$\text{Det}(\alpha^2 A + \mu^2 B + \alpha\mu C + \alpha D + \mu E + F) \quad \textbf{Q.E.D.}$$

Now we consider the bi-matrix variate, matrix quadratic equation specified in case (ii) above. It is the "dual" of the equation considered in case (i). The following Lemma is an interesting result.

**Lemma 3:** Consider the bi-matrix variate, matrix quadratic equation of the form in case (ii) above i.e.
$$AX^2 + BY^2 + CXY + DX + EY + F \equiv \bar{0}$$
Suppose the matrices $\{X, Y\}$ have common right eigenvectors corresponding to all the eigenvalues $\{\alpha, \mu\}$ (not necessarily same). Then, we necessarily have that all the eigenvalues of unknown matrices $\{X, Y\}$ are zeroes (solutions) of the determinental polynomial (bi-scalar variate)

$$\text{Det}(\alpha^2 A + \mu^2 B + \alpha\mu C + \alpha D + \mu E + F).$$

**Proof:** Using an argument similar to the one employed in Lemma 2 (by using the right eigenvectors of X, Y), the above result follows. Details are avoided for brevity. **Q.E.D.**

Suppose, we consider the case where X, Y are diagonalizable (one sufficient condition is that all the eigenvalues are distinct). If all the right eigenvectors of X, Y are same, by the argument used in proof of Lemma 1, all the left eigenvectors are same.

**Lemma 4:** Consider the bi-matrix variate, matrix quadratic equation of the form in case (iii) above i.e.
$$X A X + Y B Y + X C Y + X D + Y E + F \equiv 0$$
Suppose the matrices { X, Y } have common right eigenvectors corresponding to all the eigenvalues $\{\alpha, \mu\}$ ( not necessarily same ). Then, we necessarily have that all the eigenvalues of unknown matrices { X, Y } are zeroes ( solutions ) of the determinental polynomial ( bi-scalar variate )

$$\text{Det} ( \alpha^2 A + \mu^2 B + \alpha \mu C + \alpha D + \mu E + F ).$$

**Proof:** Using an argument similar to the one employed in Lemma 1, Lemma 2, the above result follows. Details are avoided for brevity.   Q.E.D.

**Remark 1**: Using Lemmas 2 to 4, we realize that the eigenvalues of X, Y are a subset of the zeroes of the determinental polynomial
$$\text{Det} ( \alpha^2 A + \mu^2 B + \alpha \mu C + \alpha D + \mu E + F ).$$
But, unlike in the case of uni-matrix variate, matrix quadratic equation, the bi-matrix variate, polynomial matrix need not factor into the following form:

$$( \alpha^2 A + \mu^2 B + \alpha \mu C + \alpha D + \mu E + F ) \equiv (\alpha I - X ) G ( \alpha, \mu )(\mu I - Y )$$

Thus, there could be zeroes of the determinental polynomial which are not necessarily the eigenvalues of X, Y. Hence, in the case of Bi-Matrix variate, matrix quadratic equation, we are unable to derive results that are possible in the uni-matrix variate case. The crucial factorization of polynomial matrix could hold true in some cases.

**Method for Determining the Structured Solutions X, Y:**

As in the uni-matrix variate case, we try to find the solutions of structured bi-matrix variate, matrix quadratic equation given that we know the eigenvalues of X, Y. Also X, Y are assumed to be diagonalizable. As discussed earlier, in all the three cases of interest, the left as well as right eigenvectors of solutions are identical. Thus, we have that
$$X = T F_1 T^{-1} \quad and \quad Y = T F_2 T^{-1} \ldots\ldots\ldots(2.15)$$

Now let us consider the equation as in case (ii) i.e.
$$A X^2 + B Y^2 + C X Y + D X + E Y + G \equiv \bar{0}.$$
Substituting for X, Y, we necessarily have that

$$A T F_1^2 T^{-1} + B T F_2^2 T^{-1} + C T F_1 F_2 T^{-1} + D T F_1 T^{-1} + E T F_2 T^{-1} + G \equiv \bar{0} \ldots (2.16)$$

Now extracting common $T^{-1}$ on the right hand side, we have that

$$[ A T F_1^2 + B T F_2^2 + C T F_1 F_2 + D T F_1 + E T F_2 + G T ] T^{-1} \equiv \bar{0} \ldots (2.17)$$

Since, T is non-singular, we necessarily have that

$$[A\,T\,F_1^2 + B\,T\,F_2^2 + C\,T\,F_1 F_2 + D\,T\,F_1 + E\,T\,F_2 + G\,T] \equiv \bar{0} \ldots\ldots\ldots\ldots(2.18)$$

Hence, it is clear that the two commuting, diagonalizable matrices can be determined by solving the above homogeneous system of linear equations. It should be clear that similar procedure can be utilized in the other cases of structured bi-matrix variate, matrix quadratic equations.

### 3. Multi-Matrix Variate, Matrix Polynomial Equations : Solution Techniques:

The results in the previous section naturally motivate us to consider the multi-matrix variate, matrix polynomial equations. But, once again it is clear that no reasonable results are possible if we donot restrict consideration to "structured" multi-matrix variate, matrix polynomial Equations. Thus, we are naturally led to considering the following cases:

(i) All the unknown matrix variables $\{X_1, X_2, \ldots, X_l\}$ occur on the left hand side of the known coefficient matrices. Also, all of them have common left eigenvectors corresponding to possibly different eigenvalues

(ii) All the unknown matrix variables $\{X_1, X_2, \ldots, X_l\}$ occur on the right hand side of the known coefficient matrices. Also, all of them have common right eigenvectors corresponding to possibly different eigenvalues

(iii) Some of the unknown matrix variables $\{X_1, X_2, \ldots, X_l\}$ occur on the left hand side of the known coefficient matrices and the others occur on the right hand side of the known coefficient matrices. Also, all of them have common left / right eigenvectors corresponding to possibly different eigenvalues.

Corresponding to the above three cases, we have the following multi-matrix variate, matrix polynomial equations. In these equations, the coefficients are all matrices (second order tensors) distinguished by the indices

(i)
$$\sum_{i_1=0}^{N}\sum_{i_2=0}^{N}\cdots\sum_{i_m=0}^{N} X_1^{i_1} X_2^{i_2} \ldots X_m^{i_m} A_{(i_1,i_2\ldots i_m)} + \sum_{j_1=0}^{N}\sum_{j_2=0}^{N}\cdots\sum_{j_m=0}^{N} X_1^{j_1} X_2^{j_2} \ldots X_m^{j_m} A_{(j_1,j_2\ldots j_m)} + \cdot$$
$$(i_1 + i_2 + \cdots\ldots\ldots + i_m = N) \qquad (j_1 + j + \cdots\ldots\ldots + j_m = N-1)$$

$$+ \ldots\ldots + \sum_{s=1}^{m} X_s F_s + G \equiv \bar{0} \qquad\qquad\qquad \ldots\ldots\ldots\ldots(3.1)$$

(ii)
$$\sum_{i_1=0}^{N}\sum_{i_2=0}^{N}\cdots\sum_{i_m=0}^{N} A_{(i_1,i_2\ldots i_m)} X_1^{i_1} X_2^{i_2} \ldots X_m^{i_m} + \sum_{j_1=0}^{N}\sum_{j_2=0}^{N}\cdots\sum_{j_m=0}^{N} A_{(j_1,j_2\ldots j_m)} X_1^{j_1} X_2^{j_2} \ldots X_m^{j_m} \cdot$$
$$(i_1 + i_2 + \cdots\ldots\ldots + i_m = N) \qquad (j_1 + j + \cdots\ldots\ldots + j_m = N-1)$$

$$+ \ldots\ldots + \sum_{s=1}^{m} F_s X_s + G \equiv \bar{0} \qquad\qquad\qquad \ldots\ldots\ldots\ldots(3.2)$$

(iii) In this case, some matrix unknowns occur on the left hand side of coefficient matrices whereas the others occur on the right hand side. Exact equation leads to cumbersome notation.

**Lemma 5:** Consider the multi-matrix variate, matrix polynomial equation of the form in (i) above i.e.

$$\sum_{i_1=0}^{N}\sum_{i_2=0}^{N}\cdots\sum_{i_m=0}^{N} X_1^{i_1} X_2^{i_2} \ldots X_m^{i_m} A_{(i_1,i_2\ldots i_m)} + \sum_{j_1=0}^{N}\sum_{j_2=0}^{N}\cdots\sum_{j_m=0}^{N} X_1^{j_1} X_2^{j_2} \ldots X_m^{j_m} A_{(j_1,j_2\ldots j_m)} + \cdot$$
$$(i_1 + i_2 + \cdots \ldots + i_m = N) \qquad\qquad (j_1 + j + \cdots \ldots + j_m = N-1)$$

$$+ \ldots\ldots + \sum_{s=1}^{m} X_s F_s + G \equiv \bar{0}$$

Suppose all the matrices $\{X_1, X_2, \ldots, X_m\}$ have common set of left eigenvectors corresponding to all the eigenvalues $\{\mu_1, \mu_2, \ldots, \mu_m\}$ (not necessarily same). Then, we necessarily have that all the eigenvalues of unknown matrices $\{X_1, X_2, \ldots, X_m\}$ are zeroes (solutions) of the determinental polynomial (multi-scalar variate)

$$Det\left(\sum_{i_1=0}^{N}\sum_{i_2=0}^{N}\cdots\sum_{i_m=0}^{N} \mu_1^{i_1} \mu_2^{i_2} \ldots \mu_m^{i_m} A_{(i_1,i_2\ldots i_m)} + \sum_{j_1=0}^{N}\sum_{j_2=0}^{N}\cdots\sum_{j_m=0}^{N} \mu_1^{j_1} \mu_2^{j_2} \ldots \mu_m^{j_m} A_{(j_1,j_2\ldots j_m)} + \cdot\right.$$
$$(i_1 + i_2 + \cdots \ldots + i_m = N) \qquad\qquad (j_1 + j + \cdots \ldots + j_m = N-1)$$

$$+ \ldots\ldots + \sum_{s=1}^{m} \mu_s F_s + G \bigg). \qquad\qquad\qquad \ldots\ldots\ldots\ldots(3.3)$$

**Proof:** Follows from a similar argument as in the case of Lemma 2 for the bi-matrix variate case. Q.E.D.

**Lemma 6:** Consider the multi-matrix variate, matrix polynomial equation of the form in case (ii) above i.e.

$$\sum_{i_1=0}^{N}\sum_{i_2=0}^{N}\cdots\sum_{i_m=0}^{N} A_{(i_1,i_2\ldots i_m)} X_1^{i_1} X_2^{i_2} \ldots X_m^{i_m} + \sum_{j_1=0}^{N}\sum_{j_2=0}^{N}\cdots\sum_{j_m=0}^{N} A_{(j_1,j_2\ldots j_m)} X_1^{j_1} X_2^{j_2} \ldots X_m^{j_m}.$$
$$(i_1 + i_2 + \cdots \ldots + i_m = N) \qquad\qquad (j_1 + j + \cdots \ldots + j_m = N-1)$$

$$+ \ldots\ldots + \sum_{s=1}^{m} F_s X_s + G \equiv \bar{0}$$

Suppose all the matrices $\{X_1, X_2, \ldots, X_m\}$ have common set of right eigenvectors corresponding to all the eigenvalues $\{\mu_1, \mu_2, \ldots, \mu_m\}$ (not necessarily same). Then, we necessarily have that all the eigenvalues of unknown matrices $\{X_1, X_2, \ldots, X_m\}$ are zeroes (solutions) of the determinental polynomial (multi-scalar variate)

$$Det \left( \sum_{i_1=0}^{N}\sum_{i_2=0}^{N}\cdots\sum_{i_m=0}^{N} \mu_1^{i_1} \mu_2^{i_2} \cdots \mu_m^{i_m} A_{(i_1,i_2\cdots i_m)} + \sum_{j_1=0}^{N}\sum_{j_2=0}^{N}\cdots\sum_{j_m=0}^{N} \mu_1^{j_1} \mu_2^{j_2} \cdots \mu_m^{j_m} A_{(j_1,j_2\cdots j_m)} + \right.$$

$(i_1 + i_2 + \cdots\cdots + i_m = N)$   $(j_1 + j + \cdots\cdots + j_m = N-1)$

$+ \ldots + \sum_{s=1}^{m} \mu_s F_s + G \big)$ .                                              ……………(3.4)

**Proof**: Follows the same argument as in the case of Lemma 3. Details are avoided for brevity                                                                               Q.E.D.

Similarly Lemma corresponding to case (iii) can be stated and proved.

**Remark 2:** In the case of multi-matrix variate, matrix polynomial equations, the eigenvalues of matrix variables constitute the zeroes of the associated determinental polynomial. Since factorization of the corresponding polynomial matrix doesnot necessarily hold true ( as in the uni-matrix variate case ), there can be zeroes of the determinental polynomial which are not necessarily the eigenvalues of unknown matrices.

**Remark 3:** Since all diagonalizable unknown matrices share common left and right eigenvectors, a method can easily be derived ( as in the case of bi-matrix variate, matrix quadratic equation ) for determining the unknown matrix variables.

4. **Generalizations to Tensor Linear Operator based Polynomial Equations:**

It is well known that matrices are second order tensors. Also, the concepts such as outer product, contraction, inner product of tensors are well understood. Thus using these concepts, it is possible to define and study muti-tensor variate, tensor polynomial equations. Results from tensor algebra such as eigenvalues of linear operators are well understood. It is easy to see that the results in Section 2, Section 3 can easily be generalized to tensor linear operator based equations ( using standard concepts in tensor algebra ).

5. **Conclusions:**

In this research paper, certain results on determining the eigenvalues of structured bi-matrix variate, matrix quadratic equations are proved. Given the eigenvalues, method of determining the unknown matrices is discussed. The results are generalized to structured multi-matrix variate, matrix polynomial equations. Briefly extension of results to tensor linear operator based polynomial equations is proposed.